\documentclass{article}
\usepackage{amssymb}
\usepackage{amsmath}

\usepackage[thmmarks,amsmath,standard]{ntheorem}
\usepackage{graphics}

\newcommand{\sos}{\mathbb}

\theoremsymbol{\ensuremath{_\blacksquare}}
\newtheorem{observation}{Observation}

\theoremsymbol{}
\newtheorem{conjecture}{Conjecture}

\theoremstyle{nonumberplain}

\theoremstyle{nonumberplain}
\theoremheaderfont{\scshape}
\theorembodyfont{\normalfont}
\theoremsymbol{\ensuremath{_\blacksquare}}

\newtheorem{proofOfTheorem}{Proof of Theorem \ref{thm:main}}

\title{Line-graphs of cubic graphs are normal}

\author{Zsolt Patakfalvi\\
}

\begin{document}
\maketitle

\begin{abstract}
A graph is called normal if its vertex set can be covered by cliques
$Q_1,Q_2,\dots,Q_k$ and also by stable sets $S_1,S_2,\dots,S_l$, such that
$S_i \cap Q_j \neq \emptyset$ for every $i,j$. This notion is due to
K\"orner, who introduced the class of normal graphs as an extension of the
class of perfect graphs. Normality has also relevance in information theory.
Here we prove, that the line graphs of cubic graphs are normal.
%
%
\end{abstract}

%
%
\section{Introduction}

The concept of graph normality was introduced by J\'anos K\"orner in
\cite{K73}, where he proved that all members of the celebrated class of perfect
graphs are normal. These two graph classes come up together also in the study
of additivity properties of the information theoretic functional graph
entropy, cf. \cite{CsKLMS}, \cite{KL}, \cite{Psurv}. 
Perfect graphs are important and well-studied for several
reasons (cf., e.\ g., \cite{RamReed} and \cite{Prog}), one of which is
certainly the long-standing conjecture of Claude Berge known as the Strong
Perfect Graph Conjecture, which has been proved recently by Chudnovsky,
Robertson, Seymour, and Thomas \cite{SPGT}. A conjecture of similar flavour
has been formulated by De Simone and K\"orner \cite{CdSKJ} to characterize
those graphs all induced subgraphs of which are normal. This latter conjecture is still open. (See it below.) 

There are several important subclasses of perfect graphs known, and
by K\"orner's theorem in \cite{K73} all these are also normal. Very
few graph families were, however, identified yet as being normal among those that are not necessarily perfect. The main goal of this
paper is to show that the members of a certain family, namely those graphs
that are line graphs of cubic graphs are all normal. This class of
graphs is interesting, for example, for its connections to the four colour
theorem and related conjectures, see, e.\ g., \cite{NS}.

Now we give the definition of normality. Consider a graph $G$. A set $\sos{A}$
of subsets of $V(G)$ is a covering, if every vertex of $G$ is
contained in an element of $\sos{A}$. 

\begin{definition}
\label{def:normal}
A graph $G$ is normal if there exist two coverings, $\sos{C}$ and $\sos{S}$,
of $G$, where each element of $\sos{C}$ induces a clique, each element of
$\sos{S}$ induces a stable set and $C \cap S \neq \emptyset$ for
every $C \in \sos{C}$, $S \in \sos{S}$.  
\end{definition}

From the symmetry of Definition \ref{def:normal} it follows that a graph
is normal iff its complement is normal.
If we require normality for every induced subgraph, we obtain the notion
of hereditary normality. Clearly, every perfect graph is
hereditarily normal, since every induced subgraph of a perfect graph is 
perfect and consequently normal.

The simplest graphs which are normal but not perfect are the odd cycles with
at least 9 vertices (see \cite{K73}). Smaller odd cycles are either perfect, like
the triangle, or not even normal, like the cycles with 5 or 7 vertices.
Actually, these latter graphs and the complement of the 7-cycle are the only
minimal not hereditarily normal graphs
known so far. This motivates the following conjecture formulated by De Simone
and K\"orner \cite{CdSKJ}. 

\begin{conjecture}
\label{cnj:korner}
A graph is hereditarily normal iff neither the graph nor its complement 
contains a 5-cycle or a 7-cycle as an induced subgraph.
\end{conjecture}

For partial results on Conjecture \ref{cnj:korner} see also \cite{AW}.
As De Simone and K\"orner remarks in \cite{CdSKJ}, from Conjecture
\ref{cnj:korner} it would immediately follow that the class of heredtarily normal graphs
can be recognized in polynomial time. The analogous statement for perfect
graphs is true, though very far from being trivial, cf. 
\cite{CB}, \cite{RBG}, \cite{RPG}.

%
%
%
\section{Main result}
\label{sec:result}

Recall, that
the line graph $L(G)$ of a graph $G$ is a graph whose vertices are the edges of
$G$, and two vertices of $L(G)$ are connected if and only if the edges corresponding
to them share a common vertex in $G$. Our main result is the following.

\begin{theorem}
\label{thm:main}
The line-graph of every cubic graph is normal
\end{theorem}

Let $G$ be a cubic graph whose edge-chromatic number is 3.
The colour-classes of a good edge-colouring of $G$ form stable sets
in $L(G)$.
Every 3 edges of $G$ which share a common vertex
form a clique in $L(G)$. The set of all these cliques, together with
the previously mentioned stable sets fulfill the requirements of normality.
So, Theorem \ref{thm:main} is trivial in case of such graphs.

The nontrivial part of the theorem is that the statement also holds for
the so called snarks (\cite{NS}), i. e. those cubic graphs which have edge-chromatic number 4. A similar phenomenon can be observed in case of some famous theorems and conjectures like the Four Colour Thorem, the 5-flow Conjecture or the Cycle Double Cover Conjecture. They can also be reduced to statements about snarks, and then the proof of those statments are unknown or by far non-trivial.

The paper is organized as follows. In Section 
\ref{sec:preliminaryObservations} some preliminary observations are made 
about normal graphs and normal line graphs. In Section
\ref{sec:regular2edgeConnected} we prove a technical lemma, Lemma
\ref{lma:2_connected2}, from whcih the statement of Theorem \ref{thm:main}
follows for 2-edge-connected graphs. In Section
\ref{sec:decomposition} Theorem \ref{thm:main} is proven using Lemma
\ref{lma:2_connected2} via
some decomposition statements. Finally,
in Section \ref{sec:further} further remarks are given.

\section{Preliminary observations}
\label{sec:preliminaryObservations}

In this section we make some general observations, not restricting 
ourselves to cubic graphs. All graphs are assumed to be simple, i. e.,
without loops and multiple edges.

\begin{observation}
\label{obs:clique_stable_set}
Let $G$ be a graph. The edges corresponding to the vertices of a clique of $L(G)$ form  a star or
a triangle, while the edges corresponding to a stable set of $L(G)$ form a
matching in $G$.
\end{observation}

For easier notation we introduce the following notion.

\begin{definition}
A graph $G$ is edge-normal if its line-graph $L(G)$ is normal.
\end{definition}

 In the next observation we express edge-normality, 
using Observation \ref{obs:clique_stable_set}, with the properties of $G$.
A set
$\sos{B}$ of subsets of $E(G)$ is an \emph{edge-covering}, if for every $e \in
E(G)$ there exists a $B \in \sos{B}$, for which $e \in B$.

\begin{observation}
\label{obs:edge_normal}
A graph $G$ is edge-normal iff there exist two edge-coverings, $\sos{C}$ and
$\sos{S}$, of $G$, where every member of $\sos{C}$ is a star or a
triangle, every member of $\sos{S}$ is a matching and for every $S
\in \sos{S}$ and $C \in \sos{C}$, $S \cap C \neq \emptyset$.
\end{observation}

The disturbing presence of the triangles in Observation
\ref{obs:edge_normal} motivate the introduction of our next notion.

\begin{definition} 
A graph $G$ is strongly edge-normal if there exist two edge-coverings, 
$\sos{C}$ and $\sos{S}$, of $G$, where every member of $\sos{C}$ is a 
star, every member of $\sos{S}$ is a matching and for every $S
\in \sos{S}$ and $C \in \sos{C}$, $S \cap C \neq \emptyset$.
\end{definition}

Clearly, from
strong edge-normality edge-normality follows. There are very
simple examples, which are edge-normal but not strongly edge-normal. The
simplest is the triangle itself. We will prove, that every cubic
graph is strongly edge-normal, not just edge-normal.

Now, we rewrite the definition of strong edge-normality using that a star
can be described with its midpoint. Clearly, an edge-covering $\sos{C}$
containing only stars covers the edges of $G$ iff the midpoints of the stars of
$\sos{C}$ form a vertex cover.
Recall, that $C \subseteq
V(G)$ is a \emph{vertex cover}, if for every edge $e$ of $G$ at least one endpoint of $e$
is in $C$. 
If $Y$ is a set of edges, then
we can define the \emph{endpoint set} of $Y$: $v(Y):=\{x \in V(G) : \exists y \in
V(G) , \{x,y\} \in Y \}$. If $e$ is an edge, then $v(e)=v(\{e\})$.

\begin{observation}
\label{obs:edge_normal_equiv1}
A graph $G$ is strongly edge-normal iff there exists a vertex cover 
$C$ and an edge-covering $\sos{S}$ of $G$, so that every element of
$\sos{S}$ is a matching and for every $S \in \sos{S}$, $C \subseteq v(S)$.
\end{observation}

The following equivalent version of Observation \ref{obs:edge_normal_equiv1}
will be useful.

\begin{observation}
\label{obs:edge_normal_equiv2}
A graph $G=(V,E)$ is strongly edge-normal iff there exists a vertex cover
$C \subseteq V$, so that for every $e \in E$  there is a
matching $S$, for which $C \subseteq v(S)$ and $e \in S$.
\end{observation}

From the definition of strong edge-normality it follows easily that if every component of a graph is
strongly edge-normal, then the whole graph is strongly edge-normal. Thus
from now on, we will restrict ourselves to connected graphs.

\section{Line-graphs of 2-edge-connected cubic graphs}
\label{sec:regular2edgeConnected}

In this section we prove Lemma \ref{lma:2_connected2}, from which the statement
of Theorem \ref{thm:main} follows immediately for 2-edge-connected graphs.
The proof of Theorem \ref{thm:main} will be finished
in the next section via a decomposition argument. 
There, the non-2-edge-connected cubic
graphs will be decomposed into smaller subgraphs along their cutting edges. 
At the end of the
decomposition 'almost 2-edge-connected' graphs will be obtained, the strong edge-normality of which 
we will be able to prove using the results of this section.

First some preparations are needed.
Consider a graph $G=(V,E)$ and one of its vertices $v$. $N(v):=\{x \in V 
: \{x,v\} \in E \}$ is the \emph{set of neighbours} of $v$.  
The \emph{degree} of $v$ is the number of its neighbours, $|N(v)|$.
If $U \subseteq V$, then $G[U]$ is the subgraph of $G$ induced by $U$.
A vertex cover $C$ of $G$ is \emph{good}
if it satisfies the requirements of Observation \ref{obs:edge_normal_equiv2},
that is, for every $e \in E$ there is a matching $S$, for which $C \subseteq
v(S)$ and $e \in S$. Clearly, $G$ is strongly edge-normal iff it has 
a good vertex cover.
If $C$ is a good vertex cover, then every vertex cover $C' \subseteq C$ is
also good, that is, $G$ is strongly edge normal iff it has a good
minimal vertex cover (with respect to inclusion). 
So, our first goal is to find a necessary and sufficient condition for the
existance of good minimal vertex covers.

Now, suppose that $G$ is cubic.
Consider a minimal vertex cover $C$ of $G$. If $G[C]$ contains a vertex $v$
with degree 3, then $C \setminus \{ v \}$ is also a vertex cover, that is, 
$C$ is not minimal. (Recall, that $G$
is cubic.) So, the maximum degree of $G[C]$ is 2. This explains, why 
the following notions will be useable for us in the characterization of
good minimal vertex covers.

\begin{figure}
\begin{center}
\includegraphics{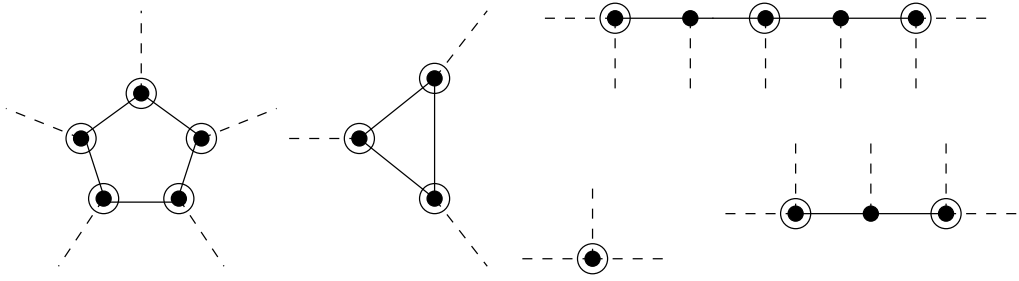}
\caption{Important vertices of odd components}
\label{pic:imp_vert}
\end{center}
\end{figure}

Consider a graph $F$ and $U \subseteq V(F)$, for which $F[U]$ has maximum
degree at most 2. Then $F[U]$ consists of three types of components: isolated points, 
cycles and paths. Consider a component $H$ of $F[U]$. 
$H$ is an \emph{odd component}, if it has an odd number of vertices.
If $H$ is an odd path with the veritces $v_1,v_2, \dots, v_n$ in their natural
order along the path, then $v_i$ is an \emph{odd vertex} of $H$ if $i$ is odd and an \emph{even}
vertex of $H$
otherwise. If $H$ is an odd component, then a vertex $v$ of $H$ is \emph{important} 
if $H$ is an isolated point or an odd cycle and $v$ is an arbitrary vertex
of $H$ or if $H$ is an odd path and $v$ is an odd vertex of $H$. 
The important vertices of odd components with at most five vertices are shown in Figure \ref{pic:imp_vert}. (The important vertices are denoted by a circle drawn around them.)

%
%
\begin{definition}
\label{def:nice_matching}
Consider a graph $F$ and $U \subseteq V(F)$ for which $F[U]$ has maximum degree
at most 2.
A matching $S \subseteq E(F)$ is nice for $U$ in $F$, if it has the following
structure.
$S=\{\{x_H,y_H\}: \textrm{$H$ is an odd component of $F[U]$}\}$,
where $x_H$ is an important vertex of $H$ and $y_H \in V(F) \setminus U$ 
for every odd component $H$ of $F[U]$.
\end{definition}

\begin{figure}
\begin{center}
\includegraphics{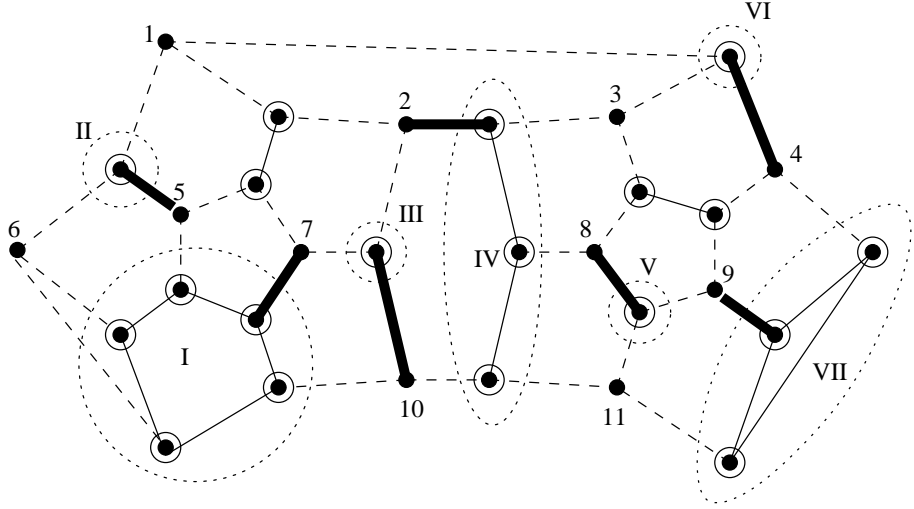}
\caption{Nice matching in a cubic graph}
\label{pic:nice_match}
\end{center}
\end{figure}

In Figure \ref{pic:nice_match} an example is shown on nice matchings, where $F$ is a cubic graph, and $U$ is a minimal vertex cover of $F$. In the figure, an extra circle is drawn around the elements of $U$, the odd components of $F[U]$ are surrounded with dotted circles and ellipses, and the edges of the nice matching are thickened. The odd components and the vertices outside of $U$ are enumerated, because later we will use them in other examples. Notice, that in compliance with Definition \ref{def:nice_matching}, in the figure to every odd component exactly one edge of the nice matching belongs, and  that edge covers an important vertex of the adequate component.

The next proposition clarifies, why the foregoing notions were introduced.
It characterizes good minimal vertex covers of cubic graphs using
nice matchings. Since a cubic graph is strongly edge-normal iff it has a
good minimal vertex cover, Proposition \ref{prop:nice} also characterizes
strongly edge-normal cubic graphs.

%
%
\begin{proposition}
\label{prop:nice}
If $C$ is a minimal vertex cover in a cubic graph $G=(V,E)$, then
$C$ is good iff for every $e=\{x,y\} \in E$ there is a nice matching for
$C \setminus \{x,y \}$ in $G[V \setminus \{x,y\}]$.
\end{proposition}

Notice, that we have already proved that $G[C]$ has maximum degree at most 2, from
which the same follows for $G[C \setminus \{x,y\}]$. That is, it does make
sense to talk about a nice matching for $C \setminus \{x,y\}$ in $G[V
\setminus \{x,y\}]$.

\begin{figure}
\begin{center}
\includegraphics{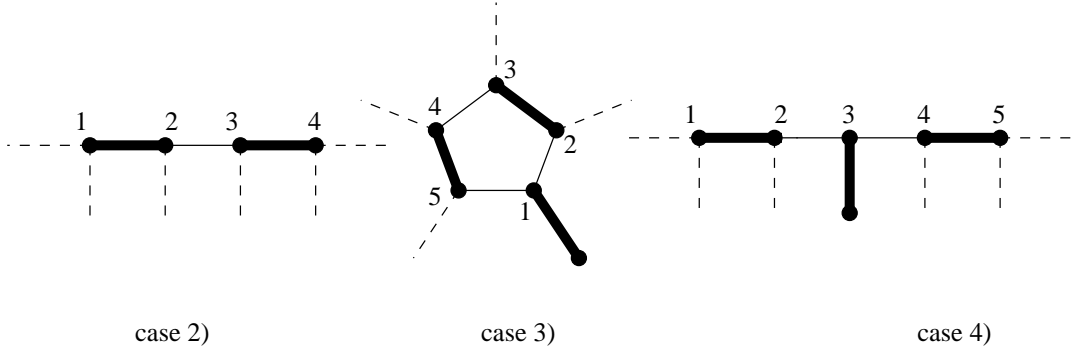}
\caption{Visualization of the proof of Proposition \ref{prop:nice}}
\label{pic:nice_match_prop}
\end{center}
\end{figure}

\begin{proof}
$\Longleftarrow$: 
We need to show, that for every edge $e$ there is a matching $S$,
for which $C \subseteq v(S)$ and $e \in S$. Fix an edge $e=\{x,y\}$.
Let $N$ be the nice matching for $C \setminus \{x, y \}$ in 
$G[C \setminus \{x,y\}]$. We
can get a matching we need, in a form of $K:=\{e\} \cup N \cup M$,
where every member of $M$ will be an edge of a component of $G[C \setminus
\{x,y\}]$. We define below for every component $H$ of $G[C \setminus
\{x,y\}]$, which edges of $H$ should be put in $M$.
There are four cases, some of which are visualized in Figure \ref{pic:nice_match_prop}. There the edges of $N$ and those, which should be put in $M$ are denoted by thick lines. The cases are the following.
%
%
\begin{enumerate}
%
%
\item $H$ is an isolated point. Since $N$ is a nice matching, $H$ is
covered by $v(N)$, so there is no need to bother with such components.
%
%
\item $H$ is an even cycle or an even path, with vertices $v_1,v_2,\dots,
v_{2k}$ in their natural order along the cycle or path. Then there
are no points of $H$ in $v(N)$, so $\{v_i,v_{i+1}\}$ should be put in
$M$ for every odd $i$, for which $1 \leq i \leq 2k-1$.
%
%
\item $H$ is an odd cycle, with vertices $v_1,v_2,\dots,v_{2k+1}$ in
circular order. Then there is one point in $H$ wich is in $v(N)$. It
can be supposed, that it is $v_1$. In this case $\{v_i,v_{i+1}\}$ should be put in $M$
for every even $i$, for which $2 \leq i \leq 2k$.
%
%
\item $H$ is an odd path, with the vertices $v_1,v_2,\dots,v_{2k+1}$
in the natural order. Then there is one point $v_j$ of $H$ in $v(N)$. From the
definition
of nice matching it follows that $j$ is odd, so the edges $\{v_i,v_{i+1}\}$ can
be put in $M$ for every $1 \leq i \leq 2k$, where $i-j$ is odd if
$i>j$, even otherwise.
\end{enumerate}

It is clear, that the edge set $K$ we obtain this way is a matching, for which $e \in K$
and $C \subseteq v(K)$.

\vspace{5pt}

$\Longrightarrow$: Fix an $e=\{x,y\} \in E$. Since $C$ is good, 
there is a matching $S$, for which $C \subseteq v(S)$ and $e \in S$. 
Clearly, $S \setminus \{e\}$ is a matching in $G[V \setminus \{x,y\}]$.
We will show, that there is a
subset of $S \setminus \{e\}$, which is a nice matching for $C \setminus
\{x,y\}$ in $G[V \setminus \{x,y\}]$.

If an edge has only one of its endpoints in a component of
$G[C \setminus \{x,y\}]$, then the other endpoint of the edge cannot be in $C \setminus
\{x,y\}$. That is, it is enough to show, that for every odd component $F$ there is an
edge $g$ of $S \setminus \{e\}$ which has only one of its endpoints $v$ in $F$, 
and if $F$ is an odd path, we can choose $g$, so that $v$ is one of the odd vertices
of $F$. Let us call such an edge $g$ a
nice edge for the component $F$. 

Now we show, that there is at least one nice edge in $S \setminus \{e\}$
for every odd component
of $G[C \setminus \{x,y\}]$. Fix such an odd component $H$. Clearly,
$V(H) \subseteq v(S \setminus \{e \})$. Since $H$ is odd,
there should be an edge in $S \setminus \{e\}$, which has only one endpoint 
in $H$. Otherwise every edge of $S \setminus \{e\}$ would have 0 or 2 of its
endpoints in $H$, that is, $H$ would have an even number of vertices.
It remains to prove only, that if $H$ is an odd
path, then there is an edge with one endpoint $v$ in $H$, such that $v$
is an odd vertex of $H$. Suppose the contrary. 
Let $v_1,v_2, \dots, v_k$ be the vertices of $H$ in the natural order.
\[
j:=\min \{ i \in \mathbb{Z} : 1 \leq i \leq k , 
\textrm{there is an edge $\{v_i,w\}$ in $S \setminus \{e\}$, so that $w \notin
V(H)$}  \}
\]
This minimum exists, since $H$ is an odd component, and it is even, since
there is no nice edge for $H$. So, 
$F:=H[\{v_1, \dots, v_{j-1}\}]$ is an odd path or an isolated point, and
it is covered with edges of $\{S \setminus e\}$, which have
both of their endpoints in $F$. This is a contradiction.
\end{proof}

%
%

Some more notations are needed for Lemma \ref{lma:2_connected2}.
Let $X$ and $Y$ be disjoint vertex sets of a graph. We say that an edge
goes from $X$ to $Y$ or connects $X$ to $Y$ if one of its endpoints is in
$X$, and the other one is in $Y$. $E(X,Y)$ is the set of edges that
connects $X$ to $Y$.
$C_k$ and $P_k$ denote the cycle and the path with $k$ vertices, respectively.
We define neighbour set and degree of vertex sets, too. Consider a graph
$G=(V,E)$ and $X \subseteq V$.
The \emph{neighbour set} of $X$ contains those vertices that are not in $X$
but have a neighbour in $X$.
The \emph{degree} of $X$ is the number of edges, which connect $X$
to $V \setminus X$. With formulae:
\[
N(X):=\{v \in V \setminus X : \exists x \in X, \{x,v\} \in E \} \quad\textrm{and} \quad 
d_G(X)=|\{\{x,y\} \in E: x \in X, y \in V \setminus X\}|
\]

The next lemma is a technical one. It states, that a graph with a 
minimal vertex cover which is not good has a very special structure.
Since every minimal
vertex cover of a non-strongly-edge-normal graph is not good, these
structural restrictions hold for non-strongly-edge-normal graphs, too.
 
\begin{lemma}
\label{lma:2_connected2}
If $G=(V,E)$ is a cubic graph and $C$ is a minimal vertex cover of $G$
which is not good, then there exists a $W \subset V$ and a
cutting edge $\{v,w\} \in E$  with the following properties.
$W$ can be partitioned into 3 subsets: 
$Z$, $Y$ and $U$, where 
$Z \subset C$, 
$U=\{x,y\} \in E$ and 
$Y \cap C = \emptyset$.
This partition can be of three different types, described below.
Furthermore, $w \in W$ and after deleting $\{w,v\}$ from $G$, $G[W]$ is one
of the remaining components.
\begin{itemize}
\item[1]{\textbf{type 1:}} $|E( U , Z )|=3$, $G[Z]$ consists of triangles and isolated
points as components and $w \in U$.
\item[2a]{\textbf{type 2a:}} $|E( U , Z )|=4$, $G[Z]$ consists of triangles and isolated
points and $w \in Y$.
\item[2b]{\textbf{type 2b:}} $|E( U , Z )|=4$, $G[Z]$ consists of triangles, isolated
points and a 3-long path and $w$ is the midpoint of the 3-long path
component of $G[Z]$.
\end{itemize}
\end{lemma}
We will refer to a set which has the properties that $W$ has in Lemma
\ref{lma:2_connected2}, as a \emph{wrong} set of $G$ with respect to $C$.
The partition $W$ has according to the statement of the lemma will be called
a $(Z,Y,U)$ partition of $W$. 

In fact, from Lemma \ref{lma:2_connected2} more than a special case of
Theorem \ref{thm:main} follows. This stronger statement is formulated in 
Corollary \ref{rmk:every_good}, and follows from the fact, that there is a
cutting edge in a cubic graph with a minimal vertex cover which is not
good.

\begin{corollary}
\label{rmk:every_good}
Every minimal vertex cover of a 2 edge-connected cubic graph is good. 
\end{corollary}

The proof of Lemma \ref{lma:2_connected2} is based on Proposition \ref{prop:nice}. 
A bipartite graph $F_e$ will be defined for an edge $e=\{x,y\}$ and $C$. 
It will turn out, that there is a matching covering a specified partite
class of $F_e$ iff there is
a nice matching for $C \setminus \{x,y\}$ in $G[C \setminus \{x,y\}]$.
That is, $e$ can be chosen so, that 
there is no such matching in $F_e$. This will cause many
restrictions on $G$, which are collected in the statement of
Lemma \ref{lma:2_connected2}.

\begin{figure}
\begin{center}
\includegraphics{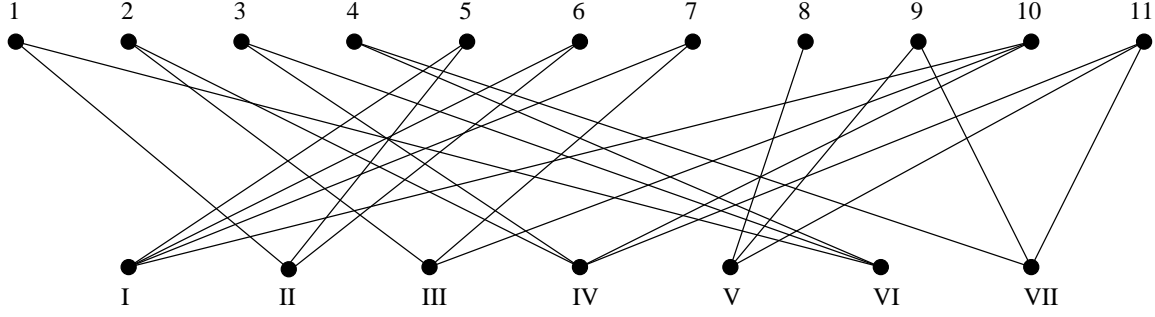}
\caption{An example on the bipartite graph of proof of Lemma \ref{lma:2_connected2}}
\label{pic:bipartite}
\end{center}
\end{figure}
\begin{proof}[of Lemma \ref{lma:2_connected2}]
Since $C$ is not good, from 
Proposition \ref{prop:nice} it follows that there must be an $e=\{x,y\} \in E$, for which
there is no nice matching for $C \setminus \{x,y\}$ in $G[V \setminus
\{x,y\}]$. Fix this edge $e$. Let $C':=C \setminus \{x,y\}$ and $V'
:= V \setminus \{x,y\}$.

Define a bipartite graph $F_e$, with partite classes $A$ and $B$, and edge set
$E_{e}$ as follows. 
\begin{itemize}
\item $A$ is the set of odd components of $G[C']$
\item $B:=V' \setminus C'$
\item $
E_{e}:=\{
\{a,b\}: \textrm{$a \in A$, $b \in B$, there is a $\{z,b\} \in E$,
where $z$ is an important  vertex of $a$.}\}.
$
\end{itemize}
For example if in the case of Figure \ref{pic:nice_match}, if we take $e$ to be the edge of a 2-long path component of the subgraph induced by the  minimal vertex cover, then $A$ consists of the odd components denoted by roman numbers on \ref{pic:nice_match}, $B$ consists of the vertices denoted by arabic numbers and $F_e$ is the graph drawn in Figure \ref{pic:bipartite}. 

Notice, that at the defnition of $F_e$ we took into consideration only those
edges on which the existence of a nice matching for $C'$ in $G[V']$ depends.
We also note, that every subset of $B$ is a stable set
in both $F_e$ and $G$. It is a stable set in $G$, because otherwise $C$
would not be a vertex cover. This property will be used later.

Clearly, there is a nice matching for $C'$ in $G[V']$ iff there is a matching
in $F_e$, which covers $A$. The selection of $e$ implies, that there is no
such matching in $F_e$. That is, the 
Hall-condition is not satisfied, so there is an $X \subseteq A$, for which
$|N(X)|<|X|$ . Fix such an $X$. Recall, that the degree of a set of vertices is
the number of edges connecting this set to the other part of the graph, and
\textrm{not} the number of vertices, to which the set is connected. Then
it is clear, that:
\begin{equation}
\label{eq:ineq_mid}
d_{F_e}(X) \leq d_{F_e}(N(X))
\end{equation}
Inequality \eqref{eq:ineq_mid} tells us, that $N(X)$ cannot be small compared to $X$.
Since $|N(X)|<|X|$, it also cannot be big. We examine how these two
conditions can be satisfied simultanously. First, we formulate necessary
conditions for \eqref{eq:ineq_mid}.
%
%
We define four useful functions: $\delta$, $\epsilon$, $t$ and $p$.
All of them are $A \rightarrow \mathbb{Z}$ functions. Consider an $a \in A$.
Then:
%
%
\[
\delta (a):= \left\{ 
  \begin{array}{ll}
  3 & \mbox{if $a$ is an isolated point of $G[C']$} \\
  2k+1 & \mbox{if $a \cong C_{2k+1}$} \\
  k+3 & \mbox{if $a \cong P_{2k+1}$}    
  \end{array}
\right.
\]
It is easy to see, that $\delta(a)$ denotes the number of edges in $G$ going from
the important vertices of $a$ somewhere outside of $a$.

%
%
\[
\epsilon (a) := 
  |\{ \{c,d\} \in E : \textrm{$c$ is an important vertex of $a$, and $d \in
  \{x,y\}$}\}| 
\]
\[
p(a) := | \{ v \in B : \textrm{there are exactly two edges of $G$,
connecting an important vertex of $a$ to $v$}\}|
\]
\[
t(a) := | \{ v \in B : \textrm{there are exactly three edges of $G$,
connecting an important vertex of $a$ to $v$}\}|
\]
For example in case of Figure \ref{pic:nice_match} and Figure \ref{pic:bipartite}, the values of $\delta$,$\epsilon$,$p$ and $t$ are displayed in the table of Figure \ref{pic:four_func} for every element of $A$.

\begin{figure}
\[
\begin{array}{r|c|c|c|c|c|c|c}
         & I & II & III & IV & V & VI & VII \\
\hline
\delta   & 5 & 3 & 3 & 4 & 3 & 3 & 3 \\
\hline
\epsilon & 0 & 0 & 0 & 0 & 0 & 0 & 0 \\
\hline
p        & 1 & 0 & 0 & 0 & 0 & 0 & 0 \\
\hline
t        & 0 & 0 & 0 & 0 & 0 & 0 & 0 \\
\end{array}
\]
\caption{Example on the values of the four functions introduced in the proof of Lemma \ref{lma:2_connected2} }
\label{pic:four_func}
\end{figure}

Consider an odd component $a \in A$ and a vertex $b \in B$. Let us call a
connection between $a$ and $b$ double or triple, if there are exactly two or three
edges, respectively, in $G$, which connect $b$ to important vertices of
$a$. We denoted by $\delta(a)$ the number of edges in $G$ going from the
important vertices of $a$ somewhere outside of $a$. We can use $\delta(a)$
to express the number of edges going out of $a$ in $F_e$, or in other words
the degree of $a$ in $F_e$, but it should be considered, that: Firstly, if there is
a double or a triple connection between $a$ and $b$, then it means two or
three edges, respectively, in $G$, but only one in $F_e$. Secondly, if there
are some edges in $G$ which go from the important vertices of $a$ to $\{x,y\}$,
they should not be
counted, because $x$ and $y$ are not in $V'$. There is no need to bother
with four time or even stronger connections, because $G$ is cubic. 
From these considerations it follows, that the degree of $a$ in $F_e$ is: 

\begin{equation}
\label{eq:dfea}
d_{F_e}(a) = \delta (a) - \epsilon (a)- p(a) - 2t(a)
\end{equation}
%
%
%
%
Since $X$ is a stable set of $F_e$, it follows from \eqref{eq:dfea} that:
\begin{equation}
\label{eq:Xdegree_first}
\begin{array}{rcl}
d_{F_e}(X) & = & \sum\limits_{a \in X} d_{F_e}(a)= \sum\limits_{a \in X} ( \delta( a) - \epsilon
(a) - p(a) - 2t(a)) = \\
& = &  \sum\limits_{a \in X} \delta(a) - \sum\limits_{a \in X} \epsilon(a) 
- \sum\limits_{a \in X} p(a) - 2\sum\limits_{a \in X} t(a) 
\end{array}
\end{equation}
We introduce easier notations for the sums in \eqref{eq:Xdegree_first}:
\[
\delta:=\sum\limits_{a \in X} \delta(a) \quad \epsilon:=\sum\limits_{a \in X}
\epsilon(a)\quad
p:=\sum\limits_{a \in X} p(a)\quad t:=\sum\limits_{a \in X} t(a) 
\]
Using our new notations, the equation obtained in \eqref{eq:Xdegree_first} can be written
as:
\begin{equation}
\label{eq:Xdegree}
d_{F_e}(X)= \delta - \epsilon - p - 2t
\end{equation}
Equation \eqref{eq:Xdegree} will be useful, by taking into acount the following facts.
\begin{enumerate}
\item $\epsilon  \leq d_G(\{x,y\})=4$
\item $\delta \geq 3|X|$, since $\delta(a) \geq 3$,
for every $a \in A$.  
\end{enumerate}
From the above it follows, that:
\begin{equation}
\label{eq:ineq_left}
d_{F_e}(X) \geq 3 |X| - 4 - p - 2t 
\end{equation}
Since $G$ is cubic, a
vertex $b$ of $B$ can have only one double or triple connection. If $b$ has
a double or a triple connection, then $d_{F_e}(b) \leq 2$ or $d_{F_e}(b)
\leq 1$, respectively. That is, if there are $p$ vertices in $N(X)$ which have double
connections, and 
$t$ vertices which have triple connections, then there are $p$
vertices, whose degree is maximum 2 in $F_e$, and $t$ vertices whose degree
is maximum 1 in $F_e$. The other vertices of $N(X)$ have degree at most 3.
Since $N(X) \subseteq B$, $N(X)$ is a stable
set in $F_e$. From the above it follows, that:
\begin{equation}
\label{eq:ineq_right}
d_{F_e} (N(X)) = \sum_{b \in N(X)} d_{F_e}(b) \leq
3 |N(X)| - p - 2t 
\end{equation}
From \eqref{eq:ineq_mid},\eqref{eq:ineq_right} and \eqref{eq:ineq_left} it
follows, that:

\begin{equation}
\label{eq:together}
3 |X| - 4 - p - 2t 
\leq d_{F_e} (X) \leq d_{F_e} (N(X)) 
\leq 3 |N(X)| - p - 2t 
\end{equation}
After subtracting the terms participating on both sides and dividing by 3:
\begin{equation}
\label{eq:xnx}
|X| - \frac{4}{3} \leq |N(X)|
\end{equation}
It is clear, that $|X|$ and $|N(X)|$ are integers, so from \eqref{eq:xnx} it
follows, that $|X| - 1 \leq |N(X)|$. Since we supposed, that $|X| > |N(X)|$:
\begin{equation}
\label{eq:x_nx}
|X| - 1 = |N(X)|
\end{equation}
Substituting \eqref{eq:x_nx} into \eqref{eq:together} we get:
\[
3 |X| - 4 - p - 2t 
\leq d_{F_e} (X) 
\leq 3 (|X|-1) - p - 2t 
\]
After some subtractions:
\[
3 |X| - 4 
\leq d_{F_e}(X) + p + 2t 
\leq 3|X| - 3
\]
Using \eqref{eq:Xdegree} we get:
\begin{equation}
\label{eq:main_ineq}
3 |X| - 4 
\leq \delta - \epsilon 
\leq 3|X| - 3
\end{equation}
We know that:
\[
0 \leq \epsilon \leq 4\ \textrm{and}\ 3|X| \leq \delta
\]
So, from \eqref{eq:main_ineq}:
\[
3 \leq \epsilon \leq 4 \quad \textrm{and} \quad 3|X| \leq \delta \leq
3|X|+1
\]
%
%
In the remaining part of the proof we talk about edges of $G$ (not
$F_e$ or anything else), so every time we count the number of edges going
from one vertex set to another, we do it in $G$.

We introduce some new notations. We denote by $Z$ the set of 
vertices of the components of $X$, and by $Z_{imp}$ the set of 
important vertices of the components of $X$. With formulae:

\[
Z:=\{z \in V: \exists a \in X, \textrm{$z$ is a vertex of $a$}\}
\]
\[
Z_{imp}:=\{z \in V: \exists a \in X, \textrm{$z$ is an
important vertex of $a$}\}
\]
Let $W:=\{x,y\} \cup N(X) \cup Z$. We will count the number of edges, that
connect $W$ to $V \setminus W$. We will obtain, that actually there is only
one such edge, which will essentially complete our proof. We will use the
following equation.
\begin{equation}
\label{eq:EWV}
|E( W , V \setminus W)| =
|E( \{x,y\} , V \setminus W)| + 
|E( N(X) , V \setminus W)| + 
|E( Z , V \setminus W)|
\end{equation}
$N(X)$ is a stable set of $G$, since every subset of $B$ is a stable set in $G$
and $N(X) \subseteq B$.
So, $d_G(N(X))=3|N(X)|$, because $G$ is cubic. That is, there
are $3|N(X)|$ edges which connect $N(X)$ to the other vertices of $G$. 
$\delta - \epsilon$ edges out of these $3|N(X)|$ edges go from $N(X)$ to $Z_{imp}$.
We know, that:
\[
 \delta - \epsilon
\geq
3|X|-4=3|N(X)|-1
\]
So, at least $3|N(X)|-1$ edges go from $N(X)$ to $Z_{imp}$, that is, $|E(N(X),
V \setminus W)| \leq 1$. $|E(N(X),V \setminus W)|=1$ can be satisfied only if
$\delta - \epsilon=3|N(X)|-1=3|X|-4$, i. e., $\epsilon=4$ and $\delta=3|X|$.

There are a few possibilities. Since we want to determine $|E(W, V \setminus
W)|$, we calculate in every case the terms appearing on the right side of
\eqref{eq:EWV}.
\begin{enumerate}
%
%
\item $\epsilon = 3$, which means according to \eqref{eq:main_ineq}, that 
$\delta = 3|X|$. So, for every component $a$ represented in $X$,
$\delta(a)=3$. From the definition of $\delta$ it follows, that $\delta(a)=3$
iff $a$ is a triangle or an isolated point.
As it was previously developed, $|E(N(X),V \setminus W)|=0$ in this case, since
$\epsilon=4$ is not satisfied. As every component
in $X$ is a triangle or an isolated point, every vertex of these components is
an important vertex, or in other words $Z=Z_{imp}$. So, every edge which has
one of its endpoints in $Z$, has the other endpoint in $N(X) \cup \{x,y\}$.
From this it follows, that $|E(Z,V \setminus W)| = 0$, too. Since $\epsilon=3$, there is an
edge, which has $x$ or $y$ as an endpoint, and the other endpoint is not in
$Z_{imp}(=Z)$. However this other endpoint neither can be in $N(X)$, because
in this case every edge which has one of its endpoints in $N(X)$ has the other
one in $Z_{imp}$, so $|E(\{x,y\},V \setminus W)|=1$.
Using \eqref{eq:EWV} we obtain $|E(W,V \setminus W)|=1$.

\begin{center}
\includegraphics{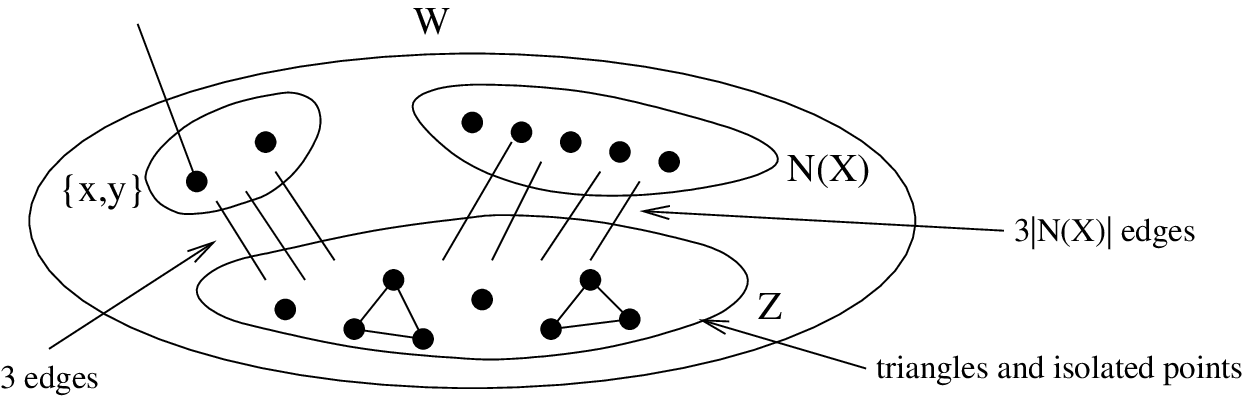}
\end{center}

%
%
\item $\epsilon = 4$, which means according to \eqref{eq:main_ineq}, that 
$\delta= 3|X|$ or $\delta=3|X|+1$, forming two subcases:
\begin{enumerate}
%
%
\item $\delta= 3|X|$, so for every component $a$ in $X$ $\delta(a)=3$, that
is, $a$ is a triangle or an isolated point.
Since $\epsilon=4$ all the four edges which has one of their endpoints in
$\{x,y\}$, and the other one outside of it, have this other endpoint in
$Z_{imp}$. So, $|E(\{x,y\},V \setminus W)|=0$. Since every component of $X$ has only important
vertices, $Z=Z_{imp}$, and $|E(Z,V \setminus W)|=0$, too.
Since all the four edges, which come out of $\{x,y\}$ go into $Z_{imp}$, the
edges which have one of their endpoint in $N(X)$ must have the other one in
$Z$ or in $V \setminus W$. Let us count how many out of these edges go from
$N(X)$ to $Z$. Since in this case $Z=Z_{imp}$, this number is:
\[
\delta - \epsilon = 3|X|-4=3|N(X)|-1
\]
So, $|E(N(X),V \setminus W)|=1$. 
Using \eqref{eq:EWV} we obtain $|E(W,V \setminus W)|=1$.

\begin{center}
\includegraphics{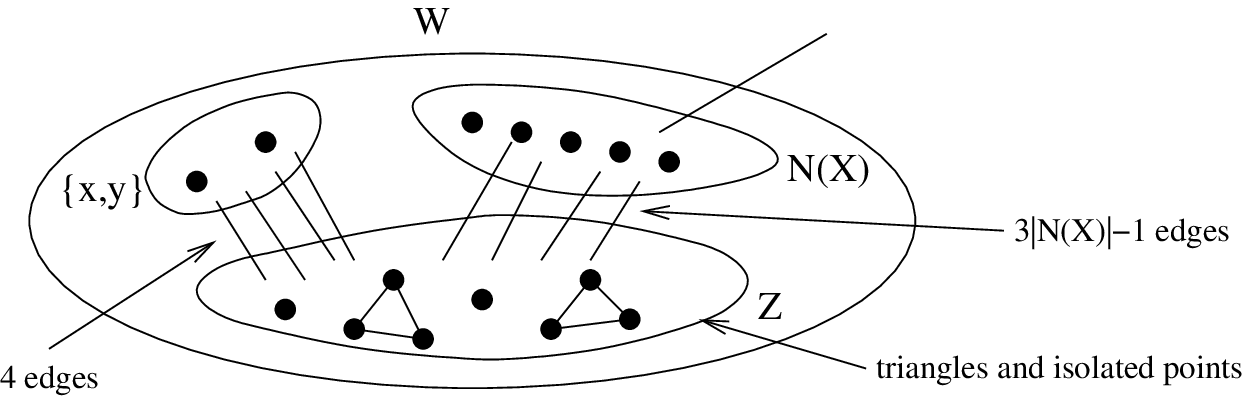}
\end{center}

%
%
\item $\delta= 3|X|+1$, so for all but one component $a$ in $X$, $\delta(a)=3$,
and there is one component $a$ for which $\delta(a)=4$. In the first case $a$
is a triangle or an isolated point, in the second case $a$ is a path of $3$
vertices.

As it was previously developped, $|E(N(X),V \setminus W)|=0$ in this case, since $\delta=3|X|$
is not satisfied. Since $\epsilon=4$, every edge from $\{x,y\}$ goes to
$Z_{imp}$, so $|E(\{x,y\},V \setminus W)|=0$. Every edge from $Z_{imp}$ goes to $N(X) \cup \{x,y\}$,
but the midpoint of the 3-long path component of $X$ is not in $Z_{imp}$. So
the edge, which goes from this vertex out of the 3-long path cannot go into
$N(X) \cup \{x,y\}$, since every edge which has one of its endpoints in $N(X)
\cup \{x,y\}$, has the other one in $Z_{imp}$, in this case. So $|E(Z,V
\setminus W)|=1$. Using \eqref{eq:EWV} we obtain $E(W,V \setminus W)=1$.

\begin{center}
\includegraphics{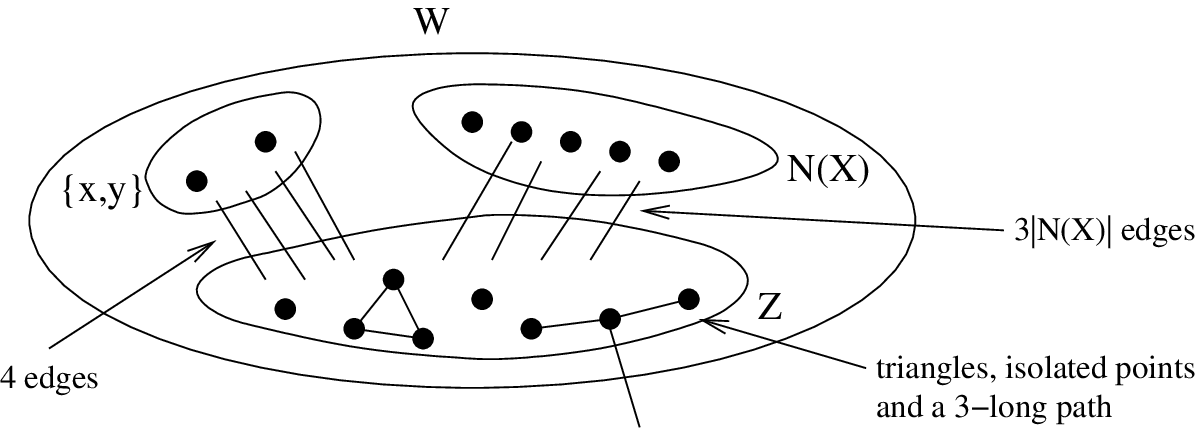}
\end{center}

\end{enumerate}
\end{enumerate}

Notice, that there are $\epsilon$ edges going from $\{x,y\}$ to $Z$ (not
just $Z_{imp}$) in every case. If we denote $N(X)$ by $Y$ and $\{x,y\}$ by
$U$, we obtain the statement of the lemma.
\end{proof}

%
%
%
\section{Decomposition}
\label{sec:decomposition}

In the previous section we have proven a special case of Theorem
\ref{thm:main}.
In this section we finish the proof of the theorem via a decomposition
argument.
We will decompose every non-strongly-edge-normal non-2-edge-connected cubic
graph along its
cutting edges into small subgraphs, such that for every graph 
at least one of the subgraphs obtained
will not be strongly-edge-normal. As a contradiction, we will prove
using our previous results, that all the small subgraphs are strongly
edge-normal.

First we define a notion which is needed because the subgraphs
obtained during the decomposition are not cubic. 

%
%
\begin{definition}
A graph $G=(V,E)$ is a cubic graph with arms, if there is a $K \subseteq V$
with the following properties.
\begin{itemize}
\item $K \neq \emptyset$ and for every $v \in K$, $d_G(v)=3$.
\item for every $v \in V \setminus K$, $d_G(v)=1$ and $N(v) \subseteq K$.
\end{itemize}
The vertices of $V \setminus K$ are the arms and $K$ is the body of
$G$.
\end{definition}

Notice, that a cubic graph is also a cubic graph with arms, because
$K=V$ is permitted in the definition. The next definition gives the
decomposition itself. Actually, we define its inverse,
the method how a graph with a cutting edge can be made by sticking together two
cubic graphs with arms. The notations of the next definition are visualized
in Figure \ref{pic:connection}.

\newcommand{\pottyx}{\circle*{7}}
\newcommand{\grafx}{\oval(70,50)}
\newcommand{\potty}[1]{\put(#1,25){\pottyx}}
\newcommand{\graf}[1]{\put(#1,25){\grafx}}

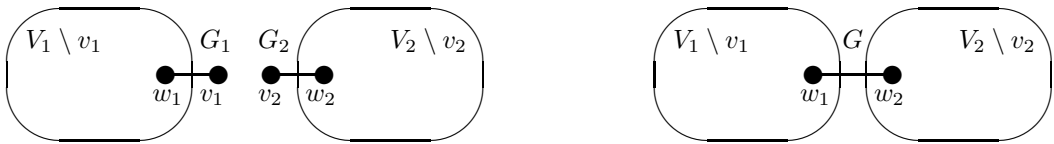
\begin{figure}[b]
\begin{center}
\begin{picture}(400,50)
\put(73,35){$G_1$}
\put(95,35){$G_2$}
\put(73,15){$v_1$}
\put(95,15){$v_2$}
\put(55,15){$w_1$}
\put(113,15){$w_2$}
\put(7,35){$V_1 \setminus v_1$}
\put(145,35){$V_2 \setminus v_2$}
\graf{35}
\potty{60}
\potty{80}
\put(60,25){\line(1,0){20}}
\put(100,25){\line(1,0){20}}
\potty{100}
\potty{120}
\graf{145}
\put(316,35){$G$}
\put(300,15){$w_1$}
\put(328,15){$w_2$}
\put(252,35){$V_1 \setminus v_1$}
\put(360,35){$V_2 \setminus v_2$}
\graf{280}
\potty{305}
\put(305,25){\line(1,0){30}}
\potty{335}
\graf{360}
\end{picture}
\caption{Two cubic graphs with arms and their connections}
\label{pic:connection}
\end{center}
\end{figure}

\begin{definition}
\label{def:connection}
Consider for $i=1$ and $2$ a cubic graph with arms $G_i=(V_i,E_i)$, one
of its arms $v_i$ and the neighbour $w_i$ of $v_i$, such that
$(V_1 \setminus \{v_1\}) \cap (V_2 \setminus \{v_2\}) = \emptyset$.
The graph $G=(V,E)$ is the connection of $G_1$ and $G_2$ at the
arms $v_1$ and $v_2$, if $V=(V_1 \setminus \{v_1\}) \cup (V_2 \setminus
\{v_2\})$, and $E=(E_1 \setminus \{\{w_1,v_1\}\})
\cup (E_2 \setminus \{\{w_2,v_2\}\}) \cup \{\{w_1,w_2\}\}$.
We denote $G$ by $G_1 +_{v_1,v_2} G_2$.
\end{definition}

If $G_1$ and $G_2$ are induced subgraphs of a graph, so that $v_1=w_2$ and
$v_2=w_1$, then $G_1 +_{v_1,v_2} G_2$ is also denoted by $G_1 +_e G_2$,
where $e$ is the edge $\{v_1,v_2\}$.
Clearly, $G_1$ and $G_2$ are part of $G$ in
a natural way. That is, there are two induced subgraphs $G_1':=G[V_1 \cup
w_2]$ and $G_2':=G[V_2 \cup w_1]$ of
$G$, isomorphic to $G_1$ and $G_2$, respectively, and $G=G_1' +_e G_2'$, 
where $e$ is the edge  $\{w_1,w_2\}$ of $G$. It is also clear, that since
$G_1$ and $G_2$ are cubic graphs with arms, $G_1
+_{v_1,v_2} G_2$ is also such a graph.

The next lemma is the heart of our decomposition. It allows us to decompose
cubic graphs with arms along their cutting edges. It states that if two
such graphs are strongly edge-normal, then their connection is also strongly
edge-normal.

%
%
\begin{lemma}
\label{lma:connection}
Let $G=(V,E)$ be a cubic graph with arms and $G_1$ and $G_2$
cubic subgraphs with arms of $G$ such that $G=G_1 +_e G_2$.
If $G_1$ and $G_2$ are strongly edge normal, then $G$ is also
strongly edge-normal.
\end{lemma}

\begin{proof} Let $U$ and $W$ be the vertex sets of the two components of
the graph obtained by deleting $\{u,w\}:=e$ from $G$. 
It can be supposed that $u \in U$, $w \in W$, $G_1=G[U \cup \{ w\}]$
and $G_2= G[W \cup \{ u\}]$. Since $G_1$ and $G_2$ are strongly edge-normal,
there is a good minimal vertex cover in both of them. Let they be
$C_1$ and $C_2$, respectively.
Clearly, for $i=1, 2$ and for every
$f \in E(G_i)$ there exists a matching $S$ in $G_i$, for which $f \in S$ and
$C_i \subseteq v(S)$. We are supposed to prove, that the same statement
holds for
$G$ and a vertex cover $C$ of $G$, instead of $G_i$ and $C_i$. 

Let $C$ be $C_1 \cup C_2$. Consider an arbitrary edge $f$ of  $G$.
$f$ is an edge of $G_1$
or $G_2$. (If it is an edge of both, we can choose.) It can be supposed, that it
is a part of $G_1$. Then there is a matching $S_1$ in $G_1$, for which $f \in
S_1$ and $C_1 \subseteq v(S_1)$.
There are two separate cases:
\begin{enumerate}
\item $\{w,u\} \in S_1$. $\{w,u\}$ is also an edge of $G_2$. Let
$S_2$ be a matching of $G_2$, for which  $\{w,u\} \in S_2$ and $C_2
\subseteq v(S_2)$. Then $S:=S_1 \cup S_2$ is a matching in $G$ we need,
because $C = C_1 \cup C_2 \subseteq v(S_1) \cup v(S_2) = v (S)$, and
$f \in S$.
%
%
\item $\{w,u\} \notin S_1$. Then there is no edge of $S_1$ in $G_2$. Let
us choose an edge $g$ of $G_2$, which has $w$ as an endpoint, but does not have $u$.
Such an edge exists, because $w$ is the neighbour of an arm of $G_2$, that is, 
$d_{G_2}(w)=3$. There
exists a matching $S_2$ of $G_2$, for which  $g \in S_2$ and $C_2 \subseteq
v(S_2)$.
$S:=S_1 \cup S_2$ is a matching in $G$, because $v(S_1) \subset U$, and 
$v(S_2) \subset W$, so $v(S_1)
\cap v(S_2) = \emptyset$. It is also true, that $C = C_1 \cup C_2 \subseteq
v(S_1) \cup v(S_2) = v (S)$ and $f \in S$. So, $S$ has the properties
required.
\end{enumerate}
\end{proof}

\begin{figure}
\begin{center}
\includegraphics{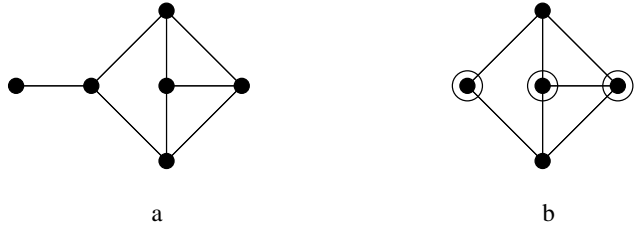}
\caption{
a) The graph used for cubic completion 
b) The cannonical vertex cover of this graph 
}
\label{pic:smallest}
\end{center}
\end{figure}

The problem with the decomposition along cutting edges is that it results
cubic graphs with arms. That is, the previously proven
Lemma \ref{lma:2_connected2} cannot be used on the results of the decomposition. 
This problem can be solved by introducing our next notion, the cubic
completion. It is a kind of smallest cubic graph which comprises a cubic
graph with arms. It arises from connecting a graph like the one in the
Figure \ref{pic:smallest}.a to the arms of a cubic graph with arms.

%
%

%
%
\begin{definition}
Consider a cubic graph with arms $G$. Let $v_1,v_2,\dots,v_k$ be
the arms of $G$, and $H_1$, $H_2$,\dots,$H_k$ graphs isomporphic to
the one in Figure \ref{pic:smallest}.a. Let $u_i$ be the only arm of $H_i$. 
The cubic completion of $G$ is the graph obtained as the connection of $G$ and 
$H_1$,$H_2$,\dots,$H_k$, if each $H_i$ is connected to $G$ at $u_i$ and $v_i$.
\end{definition}

The next lemma states that the foregoing definition is useful, that is, it is
enough to prove using Lemma \ref{lma:2_connected2} that the cubic completion 
of  a graph is strongly edge-normal. 

%
%
\begin{lemma}
\label{lma:completion}
If the cubic completion $G'$ of a cubic graph with arms $G$ is strongly edge-normal, 
then $G$ is also strongly edge-normal.
\end{lemma}

%
%
\begin{proof}
Let the arms of $G$ be $v_1,\dots,v_k$, and the neighbours of the arms be
$w_1,\dots,w_k$, respectively. $G$ is a part of $G'$ in a natural way,
so we call the common vertices of both on the names they are called in $G$.
Since $G'$ is strongly edge normal, there is a good vertex cover $C'$ in $G'$.
Then $C:=(C' \cap V(G)) \cup \{w_1,
\dots,w_k\} \setminus \{v_1,\dots,v_k\}$ is a good vertex cover in $G$, because:

Consider an arbitrary edge $e$ of $G$. There exists a matching $S$ in $G'$,
for which $e \in S$, and $C' \subset v(S)$. Let
$S':=\{s \in S: v(s) \subset V(G)\}$ and $S'':=S' \cup \{\{v_i,w_i\}:
w_i \notin v(S')\}$. Clearly, $S$ is a matching, so $S'$ and $S''$ are
also matchings. It is obvious, that $e \in S''$, and $C \subset v(S'')$
is also fulfilled.
\end{proof}
%
%

In the next lemma we will use the notions and the results of Lemma
\ref{lma:2_connected2}, which are visualized in the figures at the proof of 
Lemma \ref{lma:2_connected2}. It is clear that if a cubic graph is not normal,
there exists a wrong set of vertices, with respect to its every minimal 
vertex cover. 
So, statements which exclude the existance of a wrong set in a
cubic graph can be useful. Such statements are collected in the next lemma.

%
%
\begin{lemma}
\label{lma:technical}
Consider a cubic graph $G=(V,E)$, a  minimal vertex cover $C$, which is not
good, and a wrong set $W$ with respect to $C$. There is also a
unique edge which connects $W$ to $V \setminus W$. Let this edge be 
$\{w,v\}$, where $w \in W$. It cannot happen, that:
\begin{enumerate}
\item $w \in C$, and $N(w) \cap (C \cap W)= \emptyset$, that is, there is no
neighbour of $w$ in $C \cap W$.
\item $z,u \in W \cap C$, so that $z,u \neq w$, and $\{z,u\}$ induces a $P_2$ component in $G[C]$.
\item There is a component in $G[C \cap W]$ which is isomorphic to $P_4$.
\item There is a component in $G[C \cap W]$ which is isomorphic to $P_5$,
and one of the endpoints of this $P_5$ is $w$.
\end{enumerate}
\end{lemma}

%
%
\begin{proof}Since $W$ is a wrong set with respect to $C$, it has a
$(Z,Y,U)$ partition. Let $U=\{x,y\}$.

 First we prove, that $x \in C$ and $y \in C$ cannot be true at the
same time. Assume the contrary. Both $x$ and $y$ have two neighbours in
$V \setminus \{x,y\}$. Since $E(U,Z) \geq 3$, there is a vertex among $x$
and $y$ which has two neighbours in $Z$. It can be assumed that $x$ is this
vertex. $Z$ is a subset of $C$, so these two neighbours of $x$ are also in $C 
\setminus \{x,y\}$. That is, $x$ has all its 3 neighbours in $C$, and it is
also in $C$, which is a contradiction, since $C$ is a minimal vertex cover. 
Notice, that at least one of $x$ and $y$ must be in $C$, that is,
exactly one of $x$ and $y$ is in $C$.

It is also clear that both $x$ and $y$ must have at least one neighbour from
$Z$, since $E(U,Z) \geq 3$. If $(Z,Y,U)$ is of type 2a or 2b, then both $x$ and
$y$ has two neighbours from $Z$, since $E(U,Z) = 4$. If $(Z,Y,U)$ is of type
1, then one of $x$ and $y$ is identical to $w$, so it has one, and the other
one has two neighbours from $Z$.

The proof of the four statements follow.
\begin{enumerate}
\item Since $w \in C$, the $(Z,Y,U)$ partition can only be of type 1 or 2b.
If it is of type 1, then $w=x$ or $w=y$, so it has at least one neighbour in
$Z$, which is also in $C \cap W$. If the $(Z,Y,U)$ partition is of type 2b,
then $v$ is the second point of a 3 long path component of $G[Z]$, so it has
at least two neighbours from $C \cap W$.
%
%
\item Suppose the contrary of the statement. If $\{z,u\} \cap \{x,y\}= 
\emptyset$, then $\{z,u\}$ induces a component of $G[Z]$, but there cannot be
a $P_2$
among the components of $G[Z]$. So $\{z,u\} \cap \{x,y\} \neq \emptyset$. We
can suppose that $z \in \{x,y\}$. Then $z$ has two neighbours in $Z$, because
$z \neq w$,  so it cannot be in a component of $G[C]$, which has 2 vertices. 
%
%
\item Let us suppose the contrary of the statement. It is clear, that there
must be points from $Z$ in $P_4$ component. The components of $G[Z]$ can be
triangles, isolated points or a 3-long path. The points of triangles cannot
be part of a $P_4$, so there remain the isolated points and the 3-long path.
There is only one point apart from these components which is in $C \cap W$,
one of $x$ and $y$. We can suppose that it is $x$. It is clear, that it must
also be in the $P_4$.

If $x$ has two neighbours from $Z$, then it can connect two isolated points
or an isolated point and the 3-long path, so it is a part of a 3-long or a 
5-long path. If it has only one neighbour from $Z$, then the $(Z,Y,U)$
partition is of type 1, $x$ has only
one neighbour from $C \cap W$, which now can only be an isolated point, so
$x$ belongs to an only 2-long path component of $G[C \cap W]$, which is a
contradiction.
%
%
\item From the proof of the previous point it can be seen that a chance for
having a $P_5$ as a component of $G[C \cap W]$, is that the $(Z,Y,U)$
partition is of type 2b, $x$ (or $y$) is in $C$ and $x$ (or $y$) connects
an isolated point to the 3-long path. But in this case $w$ is the second
point of the $P_5$ (or fourth, if we start counting the vertices from the
other endpoint).
\end{enumerate}
\end{proof}
%
%

The next lemma proves, that at the end of our
decomposition we obtain graphs that are strongly edge-normal.

%
%
\begin{lemma}
\label{lma:basic}
A cubic graph with arms, which does not have a cutting edge apart from the
ones that connect its arms to its body, is strongly edge-normal.
\end{lemma}

\begin{figure}
\begin{center}
\includegraphics{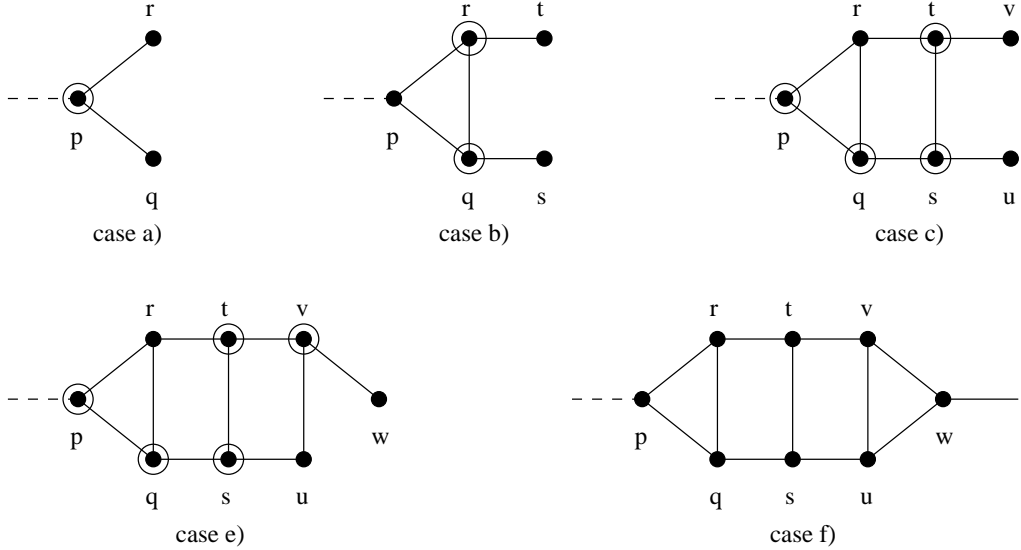}
\caption{Good vertex covers in the cubic completion of the cubic graphs with one
arm}
\label{pic:cases}
\end{center}
\end{figure}

\begin{proof}
Suppose, that there exists a non strongly edge-normal cubic graph with arms
$G$, which does
not have cutting edges apart from the ones that connect its arms to its
body. $G$ must have at least one arm, because otherwise it would be a
2-edge-connected cubic graph, which is strongly edge-normal by
Lemma \ref{lma:2_connected2}.  Let $G'=(V',E')$ be the 
cubic completion of $G$. Clearly $G'$ is also not
strongly edge-normal, because of Lemma \ref{lma:completion}.

Let $K$ be the body and $v_1,v_2,\dots,v_k$ be the arms of $G$. Let the set
of vertices by which $v_i$ was replaced in $G'$ be called $H_i$. Take a minimal
vertex cover $C$, which in each $H_i$ looks like in figure
\ref{pic:smallest}.b. ( In the figure the veritces of the vertex cover are marked with
cirles around them.) Since $G'$ is not strongly edge-normal, $C$ is not good,
so there exists a wrong  set $W$ with respect to $C$. Since $W$ is connected with only
one edge to the other parts of the graph, $W$ can only be one of the $H_i$'s, or $V'
\setminus H_i$ for some $i$. There are two cases:
\begin{enumerate}
%
%
\item $k \geq 2$. Since for every $i$ in both $G'[H_i \cap C]$ and $G'[(V'
\setminus H_i) \cap C]$,  there is a component which is isomorphic to $P_2$ and
does not contain an endpoint of a cutting edge, neither $H_i$ nor $V'
\setminus H_i$ can be a wrong set, 
because of the second statement of Lemma \ref{lma:technical}.
So, $C$ must be good for $G'$,
%
%
\item $k=1$. Since in $G'[H_1 \cap C]$ there is a component which is
isomorphic to $P_2$ and does not contain an endpoint of a cutting edge,
$H_1$ cannot be a wrong set, because of the second statement of Lemma
\ref{lma:technical}. So, $W=V' \setminus H_1$. We are going to show, that there is a minimal
vertex cover $C'$ of $G'$, for which $C' \cap H_1=C \cap H_1$, so $H_1$ is
not a wrong set with respect to $C'$, and for which neither $V' \setminus H_1$ can be a
wrong set. We will use the statements of Lemma \ref{lma:technical} to prove, that $V'
\setminus H_1$ is not a wrong set.

Let $\sos{C}$ denote the set of minimal vertex covers, which
looks like on the figure \ref{pic:smallest}.b in $H_1$. We are going to 
select from $\sos{C}$ the correct $C'$ for every type of graph. Let $p$ be
the vertex which is the endpoint of the cutting edge in $V' \setminus H_1$.
There are five cases. For each
case there is a subfigure of Figure \ref{pic:cases}, where the revelant part
of $C'$ is drawn. (In the figures, there is a circle drawn around the vertices
of $C'$.)
\begin{enumerate}
\item The two neighbours of $p$ which are in $W$ are not connected. Let us
call these two neighbours $r$ and  $q$. Then there is a $C' \in \sos{C}$, for
which $p \in C'$, and $q,r \notin C'$, because $\{q,r\}$ is a stable set in
$G'$, thus a minimal vertex cover can avoid them. This means that $V' \setminus H_1$ cannot be a wrong set, by the first
statement of Lemma \ref{lma:technical}, that is, $C'$ is good. 
\item The previous condition is not fulfilled, but the third neighbour of
$r$ (not $p$ and not $q$) is not connected to the third neighbour of $q$
(not $p$ and not $r$). Let these neighbours be $s$ and $t$, respectively.
(Though one can see easily, that $s$ and $t$ cannot coincide, it would not
matter for our reasoning.) There is a $C' \in \sos{C}$ for which $q,s,t
\notin C'$, and $q,r \in C'$, because $\{p,s,t\}$ is a stable set of $G'$. 
For this $C'$ there is a component of $G'[C' \cap W]$ which is a $P_2$ and
$p$ is not contained in this component, so $V' \setminus H_1$ cannot be a
wrong set, because of the second statement of Lemma \ref{lma:technical},
that is,  C' is good.
\item The previous conditions are not fulfilled, but the third neighbour of
$s$ (not $r$ and not $t$) is not connected to the third neighbour of $t$
(not $q$ and not $s$). Let these neighbours be $u$ and $v$, respectively.
(Here it also does not matter if $u$ and $v$ coincide.) There is a $C' \in
\sos{C}$ for which $r,u,v \notin C'$, and $p,q,s,t \in C'$, because
$\{r,u,v\}$ is a stable set in $G'$. But for this $C'$ there is a component
of $G'[C' \cap W]$ which is a $P_4$, so $V' \setminus H_1$ cannot be a wrong
set, because of the third statement of Lemma \ref{lma:technical}, that is,  C' is good.
%
%
\item The previous conditions are not fulfilled, but the third neighbour of
$u$ (not $s$ and not $v$) is not the same as the third neighbour of $v$ (not
$t$ and not $u$). Let this neighbour of $v$ be $w$. Then there is a $C' \in
\sos{C}$ for which $r,w,v \notin C'$, and $p,q,s,t,v \in C'$, because
$\{r,t,w\}$ is a stable set of $G'$. But for this
$C'$ there is a component of $G'[C' \cap W]$ which is a $P_5$, and $p$ is one
of its endpoint, so $V' \setminus H_1$ cannot be a wrong set, by the fourth
statement of Lemma \ref{lma:technical}, that is,  C' is good.
%
%
\item The previous conditions are not fulfilled. This cannot happen, because the
third edge going out of $w$ is also a cutting edge, which contradicts the
fact that $V' \setminus H_1$ is the body of $G$, and this body is 2-edge-connected.
\end{enumerate}
\end{enumerate}

So, $G'$ is strongly edge-normal, which is a contradiction.
\end{proof}
%
%

We reached the point, where we can complete the proof of Theorem \ref{thm:main}. We
make the previously promised decomposition in this proof.

%
%
\begin{proofOfTheorem}
We prove the stronger statement: Every cubic graph is strongly
edge-normal. Suppose, that there exists a graph $G=(V,E)$, which
contradicts this. Then $G$ is a non-strongly-edge-normal cubic graph
with arms. That is, there exists a minimal, nonempty subset $U$ of $V$,
so that, $G[U]$ is a non-strongly-edge-normal cubic graph with arms.
$G[U]$ does not have
cutting edges apart from the ones that connect its arms to its body, because
of Lemma \ref{lma:connection}. So, from Lemma \ref{lma:basic} it follows, that
$G[U]$ is strongly edge-normal, which is a contradiction.
\end{proofOfTheorem}

\section{Further remarks}
\label{sec:further}

First a corollary, which is the immediate consequence of the previously proven
lemmas.

\begin{corollary}
Every graph which has vertices of degree 3 or 1 only is strongly
edge-normal, that is, its line graph is normal.
\end{corollary}

\begin{proof}
Take a graph $G=(V,E)$, the vertices of which have degrees 3 and 1 only.
Remember, that it is enough to consider only connected graphs. So, if
there are no vertices with degree 3, then $G$ is a $P_2$, which is strongly
edge-normal. If there are also vertices with degree 3, then $G$ is a
cubic graph with arms. Take its cubic completion $G'$. From Theorem
\ref{thm:main} it
follows, that $G'$ is strongly edge-normal. So, $G$ is also strongly edge
normal because of Lemma \ref{lma:completion}, which proves the corollary.
\end{proof}

A graph $G$ is strongly edge-normal iff it has a good vertex
cover $C$. Since every cubic graph is strongly edge-normal, every
cubic graph has a good vertex cover. The proof we presented
is constructive in the following sense: it gives an algorithm for finding this good
vertex cover.

\begin{corollary}
There is a polynomial time algorithm, which finds a good vertex cover
for every cubic graph.
\end{corollary}

\begin{proof}
We describe the algorithm: Consider a cubic graph $G=(V,E)$. We
examine whether there exists a cutting edge in $G$. If there is
no cutting edge, then $G$ is 2-edge-connected. In this case,  as it was
developed earlier in Corollary \ref{rmk:every_good}, every minimal vertex cover
is good. Thus, the algorithm has to do only one thing, it has to find a minimal
(with respect to inclusion) vertex cover. This together with the search for
cutting edge can clearly be done in polynomial time.

From now on, suppose that $G$ is not 2-edge-connected. We decompose $G$ 
recursively along its cutting edges. We obtain induced cubic 
subgraphs of $G$ with arms
$G_1,G_2,\dots,G_k$, so that, $G$ is the connection of these subgraphs
and every subgraph has a 2-edge-connected body. 
Every $G_i$ has a body with at least one vertex, and
every vertex of $G$ is part of exactly one body, so we cannot
get more than $|V|$ graphs at the end of the decomposition. Thus, the
decomposition is finished in polynomial time. In Lemma \ref{lma:connection}
it is described how a good vertex cover of $G$ can be obtained if a
good vertex cover of every $G_i$ is given. That is, if we can find a good vertex
cover for a cubic graph with arms the body of which is 2-edge-connected
in polynomial time then we are ready. Let us fix a graph $H$ like this.

The proof of Lemma \ref{lma:basic} is constructive, that is, it distinguishes
between two
cases and four more subcases of the second case, and tells us, how a good vertex
cover of the cubic completion of $H$ should look like in each case and subcase. One
can decide in polynomial time, which case a graph $H$ belongs to, and
then finding a good vertex cover of the cubic completion is trivial. From Lemma \ref{lma:completion}
we know how to get a good vertex cover of $H$ from a good vertex cover of the
cubic completion of $H$. This can also be made in polynomial time. So the statement
of the corollary is proven.
\end{proof}

\end{document}